\documentclass[11pt]{article}

\usepackage[cp850]{inputenc}
\usepackage{latexsym,graphicx}

\font\tenmath=msbm10 scaled 1200
\font\sevenmath=msbm7 scaled 1200
\font\fivemath=msbm5 scaled 1200

\newfam\mathfam \textfont\mathfam=\tenmath
\scriptfont\mathfam=\sevenmath \scriptscriptfont\mathfam=\fivemath

\def\R{{\it I\hspace{-0.12cm}R}}

\def\N{{\it I\hspace{-0.12cm}N}}
\def\E{{\it I\hspace{-0.12cm}E}}

\def\P{{\it I\hspace{-0.12cm}P}}
\def\Z{{\it Z\hspace{-0.16cm}Z}}

\def \^#1{\if#1i{\accent"5E\i}\else{\accent"5E#1}\fi}

\def \ms{\medskip}
\def \ss{\smallskip}

\newtheorem{Thm}{Theorem}

\newtheorem{Lem}{Lemma}
\newtheorem{Pro}{Proposition}
\newtheorem{Cor}{Corollary}

\author{\sc
{\sc Harald Luschgy}\thanks{Universit\"at Trier, FB IV-Mathematik, D-54286 Trier, Germany.
E-mail: {\tt luschgy@uni-trier.de}} \quad {\sc  and}
\quad {\sc Gilles Pag\`es} \thanks{Laboratoire de Probabilit\'es et Mod\`eles al\'eatoires, UMR~7599, Universit\'e Paris 6, case 188, 4,
pl. Jussieu, F-75252 Paris Cedex 5. E-mail:{\tt  gpa@ccr.jussieu.fr}}
}
\date{February 2009}
\title{\bf Expansions for Gaussian processes and Parseval frames}

\begin{document}

%\begin{center}
%{{\bf \Large Quantization of probability distributions under } } \\ $ $ \\
%{{\bf \Large norm-based distortion measures}} \\
%\end{center}

\maketitle
\begin{abstract}

We derive a precise link between series expansions of Gaussian random vectors in a Banach space and Parseval frames in their reproducing kernel Hilbert space. The results are applied to pathwise continuous Gaussian processes and a new optimal expansion for fractional Ornstein-Uhlenbeck processes is derived. In the end an extension of this result to Gaussian stationary processes with convex covariance function is established.

\end{abstract}

\bigskip
\noindent {\em Key words:
Gaussian process, series expansion, Parseval frame, optimal expansion, fractional Ornstein-Uhlenbeck process.  }

\bigskip
%\ni {\em }

\section{Introduction}
\setcounter{equation}{0}
\setcounter{Assumption}{0}
\setcounter{Theorem}{0}
\setcounter{Proposition}{0}
\setcounter{Corollary}{0}
\setcounter{Lemma}{0}
\setcounter{Definition}{0}
\setcounter{Remark}{0}
Series expansions is a classical issue in the theory of Gaussian measures (see [2], [9], [17]).
Our motivation for a new look on this issue finds its origin in recent new expansions for fractional Brownian motions
(see [14], [1], [5], [6], [7]).

Let $(E, ||\cdot ||)$ be a real separable Banach space and let $X : (\Omega, A, \P) \rightarrow E$ be a centered Gaussian random vector with distribution $\P_X$. In this article we are interested in series expansions of $X$ of the following type. Let
$\xi_1, \xi_2, \ldots$ be i.i.d. $N(0,1)$-distributed real random variables.
A sequence $(f_j)_{j \geq 1} \in E^\N$ is called {\em admissible} for $X$ if
\begin{equation}
\sum\limits^\infty_{j = 1} \xi_j f_j \; \mbox{converges a.s. in} \; E
\end{equation}
and
\begin{equation}
X \stackrel{d}{=} \sum\limits^\infty_{j = 1} \xi_j f_j .
\end{equation}
By adding zeros finite sequences in $E$ may be turned into infinite sequences and thus also serve as admissible sequences.

We observe a precise link to frames in Hilbert spaces. A sequence $(f_j)_{j \geq 1}$ in a real separable Hilbert space
$(H, ( \cdot, \cdot))$ is called {\em Parseval frame} for $H$ if
$\sum\limits^\infty_{j = 1} (f_j, h) f_j$ converges in $H$ and
\begin{equation}
\sum\limits^\infty_{j = 1} (f_j, h) f_j = h
\end{equation}
for every $h \in H$.
Again by adding zeros, finite sequences in $H$ may also serve as frames. For the background on frames the reader is referred to [4]. (Parseval frames correspond to tight frames with frame bounds equal to 1 in [4].)
\begin{Thm}
Let $(f_j)_{j \geq 1} \in E^\N$. Then $(f_j)$ is admissible for $X$ if and only if $(f_j)$ is a Parseval frame for the reproducing kernel Hilbert space of $X$.
\end{Thm}

We thus demonstrate that the right notion of a ''basis'' in connection with expansions of $X$ is a Parseval frame and not an orthonormal basis for the reproducing kernel Hilbert space of $X$.
The first notion provides the possibility of redundancy and is more flexible as can be seen e.g. from wavelet frames. It also reflects the fact that ''sums'' of two (or more) suitable scaled expansions of $X$ yield an expansion of $X$.

The paper is organized as follows. In Section 2 we investigate the general Banach space setting
in the light of frame theory and provide the proof of Theorem 1. Section 3 contains applications to pathwise continuous processes
$X = (X_t)_{t \in I}$ viewed as $C(I)$-valued random vectors where $I$ is a compact metric space.
Furthermore, we comment on optimal expansions. Fractional Brownian motions serve as illustration.
Section 4 contains a new optimal expansion for fractional Ornstein-Uhlenbeck processes.

It is convenient to use the symbols $\sim$ and $\approx $ where $a_n \sim b_n$ means $a_n/b_n \rightarrow 1$ and $a_n \approx b_n$ means $0 < \liminf a_n/b_n \leq \limsup a_n/b_n < \infty$.
\section{The Banach space setting}
\setcounter{equation}{0}
\setcounter{Assumption}{0}
\setcounter{Theorem}{0}
\setcounter{Proposition}{0}
\setcounter{Corollary}{0}
\setcounter{Lemma}{0}
\setcounter{Definition}{0}
\setcounter{Remark}{0}
Let $(E, \|  \cdot \| )$ be a real separable Banach space. For $u \in E^*$ and $x \in E$, it is convenient to write
\[
\langle  u, x \rangle
\]
in place of $u(x)$. Let $X : ( \Omega, {\cal A}, \P) \rightarrow E$ be a centered Gaussian random vector with distribution $\P_X$. The covariance operator $C = C_X$ of $X$ is defined by
\begin{equation}
C : E^* \rightarrow E, \; C u := \E \langle u, X \rangle X .
\end{equation}
This operator is linear and (norm-)continuous. Let $H = H_X$ denote the reproducing kernel Hilbert space (Cameron Martin space) of the symmetric nonnegative definite kernel
\[
E^* \times E^* \rightarrow \R, (u, v) \mapsto \langle u, C v \rangle
\]
(see [17], Propositions III.1.6. and III.1.7).
Then $H$ is a Hilbert subspace of $E$, that is $H \subset E$ and the inclusion map  is continuous. The reproducing property reads
\begin{equation}
(h, C u)_H = \langle u, h \rangle, u \in E^*, h \in H
\end{equation}
where $(\cdot, \cdot)_H$ denotes the scalar product on $H$ and the corresponding norm is given by
\begin{equation}
\| h \|_H = \sup \{ \mid \langle u, h \rangle \mid  : u \in E^* , \langle u, C  u \rangle \leq 1 \} .
\end{equation}
In particular, for $h \in H$,
\begin{equation}
\| h \|  \leq \sup_{\| u \|  \leq 1} \langle u, C u \rangle^{1/2} \| h \|_H = \| C \|^{1/2} \| h \|_H .
\end{equation}
The $\| \cdot \|_{H^{-}}$closure of $A \subset H$ is denoted by $\overline{A}^{(H)}$.
Furthermore, $H$ is separable, $C(E^*)$ is dense in $(H, \|  \cdot \|_H)$, the unit ball
\[
U_H := \{ h \in H : \| h \|_H \leq 1 \}
\]
of $H$ is a compact subset of $E$,
\[
\mbox{supp} (\P_X) = (\mbox{ker} C)^\bot :=  \{ x \in E : \; \langle u, x \rangle= 0 \; \mbox{for every} \; u \in \; \mbox{ker} \; C \} = \overline{H} \; \mbox{in} \; E
\]
and
\begin{equation}
H = \{ x \in E : \; \|  x \|_H < \infty \}
\end{equation}
where $||x||_H$ is formally defined by (2.3) for every $x \in E$.
As for the latter fact, it is clear that $||h||_H < \infty$ for $h \in H$. Conversely, let $x \in E$ with
$||x||_H < \infty$. Observe first that $x \in \overline{H}$. Otherwise, by the Hahn-Banach theorem, there exists $u \in E^{*}$ such that $u| \overline{H} = 0$ and $\langle u, x \rangle > 0$.
Since, $\langle u, C u \rangle = 0$ this yields
\[
|| x ||_H \geq \sup_{a > 0} a \langle  u, x \rangle = \infty,
\]
a contradiction. Now consider $C(E^{*})$ as a subspace of $(H, || \cdot ||_H)$ and define
$\varphi : C(E^{*}) \rightarrow \R$ by $\varphi(C u) := \langle u,x \rangle$. If $C u_1 = C u_2$, then using
(ker$C)^{\perp } = \overline{H}, \langle u_1 - u_2, x \rangle = 0$. Therefore, $\varphi$ is well defined. The map $\varphi$ is obviously linear and it is bounded since
\[
|| \varphi || = \sup \{ | \varphi (C u) | : u \in E^{*}, ||C u ||_H \leq 1 \} = || x ||_H < \infty
\]
by (2.2). By the Hahn-Banach theorem there exists a linear bounded extension
$\tilde{\varphi} : \overline{C(E^{*})} ^{(H)} \rightarrow \R$ of $\varphi$. Then, since
$\overline{C(E^{*})} ^{(H)} = H$, by the Riesz theorem there exists $g \in H$ such that
$\tilde{\varphi} (h) = (h, g)_H$ for every $h \in H$. Consequently, using (2.2),
\[
\langle u, x \rangle = \varphi(C u) = (C u, g)_H = \langle u, g \rangle
\]
for every $u \in E^{*}$ which gives $x = g \in H$. \\

The key is the following characterization of admissibility. It relies on the Ito-Nisio theorem.
Condition $(v)$ is an abstract version of Mercer's theorem (cf. [15]. p. 43).
Recall that a subset $G \subset E^{*}$ is said to be separating if for every $x, y \in E, x \not= y$ there exists $u \in G$ such that $\langle u, x \rangle \not= \langle u, y \rangle$.
\begin{Lem} Let $(f_j)_{j \geq 1} \in E^\N$. The following assertions are equivalent.

\ss
\noindent $(i)$ The sequence $(f_j)_{j \geq 1}$ is admissible for $X$.

\ms
\noindent $(ii)$
There is a separating linear subspace $G$ of $E^*$ such that for every $u \in G$,
\[
( \langle u, f_j \rangle )_{j \geq 1 } \; \mbox{is admissible for} \; \langle u,X \rangle.
\]

\ss
\noindent $(iii)$
There is a separating linear subspace $G$ of $E^*$ such that for every $u \in G$,
\[
\sum\limits^{\infty}_{j=1} \langle u, f_j\rangle^2 = \langle u, C u \rangle .
\]

\noindent $(iv)$ For every $u \in E^*$,
\[
\sum\limits^\infty_{j=1} \langle u, f_j \rangle f_j = C u .
\]

\noindent$(v)$
For every $a > 0$,
\[
\sum\limits^\infty_{j=1} \langle u, f_j \rangle \langle v, f_j \rangle = \langle u, C v \rangle
\]
uniformly in $u, v \in \{ y \in E^{*} : \| y \|  \leq a \}$ .
\end{Lem}
{\bf Proof.} Set $X_n := \sum\limits^n_{j=1} \xi_j f_j$. (i) $\Rightarrow$ (v). $X_n$ converges a.s. in $E$ to some $E$-valued random vector $Y$, say, with $X \stackrel{d}{=} Y$. It is well known that this implies $X_n \rightarrow Y$ in $L^2_E$. Therefore,
\[
\mid \sum\limits^n_{j=1} \langle u, f_j \rangle \langle v, f_j \rangle  -  \langle u, C v \rangle \mid =
\mid \E \langle u, X_n \rangle \langle v, X_n \rangle - \E \langle u, Y \rangle \langle v, Y \rangle \mid
\]
\[
\quad \quad \quad \quad =  \mid  \E \langle u, Y - X_n \rangle \langle v, Y - X_n \rangle \mid \leq a2 \E \| Y - X_n
\|^2 \rightarrow 0 \; \mbox{as} \; n \rightarrow \infty
\]
uniformly in $u,v \in \{ y \in E^{*} : \| y \|  \leq a \}$. (v) $\Rightarrow$
(iv) $\Rightarrow$ (iii) is obvious.
(iii) $\Rightarrow$ (i). For every $u \in G$,
\[
\E \exp (i \langle u, X_n \rangle )  =  \exp ( - \sum\limits^n_{j=1} \langle u, f_j \rangle^2 /2 ) \rightarrow \exp (-\langle u, C u \rangle / 2 ) \\
 =  \E \exp (i\langle u, X \rangle).
\]
The assertion (i) follows from the Ito-Nisio theorem (cf. [17], p. 271).
(i) $\Rightarrow$ (ii) $\Rightarrow$ (iii) is obvious. $\Box$\\

Note that the preceding lemma shows in particular that $(f_j)_{j \geq 1}$ is admissible for $X$
if and only if $(f_{\sigma(j)})_{j \geq 1}$ is admissible for $X$
for (some) every permutation $\sigma$ of $\N$ so that $\sum \xi_j f_j$ converges unconditionally a.s. in $E$ for such sequences and all the $a.s.$ limits under permuations of $\N$ have distribution $\P_X$.

It is also an immediate consequence of Lemma 1(v) that admissible sequences $(f_j)$ satisfy
$\| f_j \| \rightarrow 0$ since by the Cauchy criterion,
$\lim_{j \to \infty} \sup_{||u|| \leq 1} <\langle u, f_j\rangle^2 = 0$. \\

The corresponding lemma for Parseval frames reads as follows.
\begin{Lem}
Let $(f_j)_{j \geq 1}$ be a sequence in a real separable Hilbert space $(K, (\cdot, \cdot)_K)$.
The following assertions are equivalent.
\begin{itemize}
\item[(i)]
The sequence $(f_j)$ is a Parseval frame for $K$.
\item[(ii)]
For every $k \in K$,
\[
\lim_{n \to \infty} || \sum\limits^n_{j=1} (k, f_j)_K f_j ||_K = || k ||_K .
\]
\item[(iii)]
There is a dense subset $G$ of $K$ such that for every $k \in G$,
\[
\sum\limits^\infty_{j = 1} (k , f_j)^2_K = || k ||^2_K .
\]
\item[(iv)]
For every $k \in K$,
\[
\sum\limits^\infty_{j = 1} (k, f_j)^2_K = || k ||^2_K .
\]
\end{itemize}
\end{Lem}
{\bf Proof.}
(i) $\Rightarrow$ (ii) is obvious. (ii) $\Rightarrow$ (iv).
For every $k \in K, n \in \N$,
\[
\begin{array}{lcl}
0 & \leq & || \sum\limits^n_{j=1} (k, f_j)_K f_j - k ||^2_K \\
    & = &  || \sum\limits^n_{j=1} (k, f_j)_K f_j ||^2_K - 2 \sum\limits^n_{j=1} (k, f_j)^2_K + || k ||^2_K
\end{array}
\]
so that
\[
2 \sum\limits^n_{j=1} (k, f_j)^2_K \leq || \sum\limits^n_{j=1} (k, f_j)_K f_j ||^2_K + || k ||^2_K .
\]
Hence
\[
\sum\limits^\infty_{j=1} (k, f_j)^2_K \leq ||k||^2_K .
\]
Using this inequality we obtain conversely for $k \in K, n \in \N$
\[
\begin{array}{lcl}
|| \sum\limits^n_{j=1} (k, f_j)_K f_j ||^2_K
& = & \sup_{||g||_K \leq 1} ( g, \sum\limits^n_{j=1} (k, f_j)_K f_j)^2_K \\
& = & \sup_{||g||_K \leq 1 } ( \sum\limits^n_{j=1} (k, f_j)_K (g, f_j)_K)2 \\
& \leq & \sum\limits^n_{j=1} (k, f_j)^2_K \sup_{||g||_K \leq 1} \sum\limits^n_{j=1} (g, f_j)^2_K \\
& \leq & \sum\limits^n_{j=1} (k, f_j)^2_K .
\end{array}
\]
Hence
\[
||k||^2_K \leq \sum\limits^\infty_{j=1} (k, f_j)^2_K .
\]
(iv) $\Rightarrow$ (iii) is obvious. (iii) $\Rightarrow$ (i).
Since $G$ is dense in $K$, for $k \in K$ there exist $k_n \in G$ satisfying $k_n \rightarrow k$ so that
$\lim_{n \to \infty} (k_n , f_j)^2_K = (k, f_j)^2_K$ for every $j$. Fatou's lemma for the counting measure in $\N$ implies
\[
\begin{array}{lcl}
\sum\limits^\infty_{j=1} (k, f_j)^2_K
& \leq & \liminf_{n \to \infty} \sum\limits^\infty_{j=1} (k_n, f_j)^2_K \\
  & = &  \lim_{n \to \infty} || k_n ||^2_K = ||k||^2_K .
\end{array}
\]
Therefore, one easily checks that $\sum\limits^\infty_{j = 1} c_j f_j$ converges in $K$ for every
$c = (c_j) \in l_2(\N)$ and
\[
T : l_2(\N) \rightarrow K, T(c) := \sum\limits^\infty_{j = 1} c_j f_j
\]
is linear and continuous (see [4], Theorem 3.2.3).
Consequently, the frame operator
\[
TT^{*} : K \rightarrow K, TT^{*} k = \sum\limits^\infty_{j = 1} (k, f_j)_K f_j
\]
is linear and continuous. By (ii),
\[
(TT^{*} k, k)_K = \sum\limits^\infty_{j = 1} (k, f_j)^2_K = || k ||^2_K
\]
for every $k \in G$ and thus $(TT^{*} k, k)_K = ||k||^2_K$ for every $k \in K$. This implies $TT^{*} k = k$ for every $k \in K$. \hfill{$\Box$} \\

The preceeding lemma shows that the series (1.3) converges unconditionally. Note further that a Parseval frame $(f_j)$ for $K$ satisfies
$\{ f_j : j \geq 1 \} \subset U_K$,
since
\[
||f_m||^4_K + \sum\limits_{j \not= m} (f_m, f_j)^2_K = \sum\limits^\infty_{j=1} (f_m, f_j)^2_K = || f_m ||^2_K,
\]
$\overline{span} \{ f_j : j \geq 1\} = K$
and it is an orthonormal basis for $K$ if and only if $||f_j||_K = 1$ for every $j$. \\ \\
{\bf Proof of Theorem 1.} The ''if'' part is an immediate consequence of the reproducing property (2.2) and
Lemmas 1 and 2 since for $u \in E^{*}$,
\[
\sum\limits^\infty_{j = 1} \langle u, f_j \rangle^2 = \sum\limits^\infty_{j = 1} (C u, f_j )^2_H = || C u ||^2_H = \langle u, C u \rangle.
\]
The ''only if'' part. By Lemma 1,
\[
||f_j ||_H = \sup\{ | \langle u, f_j \rangle | : \langle u, C u \rangle \leq 1 \} \leq 1
\]
so that by (2.5), $\{ f_j : j \geq 1 \} \subset H$. Again the assertion follows immediately from (2.2) and
Lemmas 1 and 2 since $C(E^{*})$ is dense in $H$. \hfill{$\Box$} \\

The covariance operator admits factorizations $C = SS^*$, where $S : K \rightarrow E$ is a linear continuous operator and
$(K, (\cdot, \cdot)_K)$ a real separable Hilbert space, which provide a useful tool for expansions. It is convenient to allow that $S$ is not injective. One gets
\begin{eqnarray}
S(K) & = & H , \\
*[.4em] (S k_1, Sk_2)_H & = & (k_1, k_2)_K , k_1 \in K, k_2 \in (\mbox{ker} S)^\bot  , \nonumber \\
*[.4em]  \| S \| & = & \| S^* \|  = \| C \|^{1/2} , \nonumber \\
*[.4em]  \overline{S^*(E^*)} & = & ( \mbox{ker} S)^\bot \; \mbox{in} \; K , \nonumber \\
*[.4em]  (\mbox{ker} S^*)^\bot & := & \{x \in E : \langle u, x \rangle = 0 \; \forall u \in \mbox{ker} S^*
\} = \overline{H} \; \mbox{in} \; E. \nonumber
\end{eqnarray}
Notice that factorizations of $C$ correspond to linear continuous operators $T : K \rightarrow H$
satisfying $TT^{*} = I$ via $S = JT$, where $J : H \rightarrow E$ denotes the inclusion map.

A sequence $(e_j)$ in $K$ is called {\em Parseval frame sequence} if it is a Parseval frame for
$\overline{span} \{e_j : j \geq 1\}$.
\begin{Pro} Let $C = S S^*, S : K \rightarrow E$ be a factorization of $C$ and let $(e_j)$
be a Parseval frame sequence in $K$ satisfying $(\mbox{ker} S)^\bot \subset \; \overline{span} \{ e_j : j = 1,2, \ldots\}$. Then $(S(e_j))$ is admissible for $X$.
Conversely, if $(f_j)$ is admissible for $X$ then there exists a Parseval frame sequence $(e_j)$ in $K$ satisfying $(\mbox{ker} S)^\perp = \overline{\mbox{span}} \{ e_j : j = 1 , 2, \ldots \}$ such that
$S(e_j) = f_j$ for every $j$.
\end{Pro}
{\bf Proof.} Let $K_0 : = \overline{\mbox{span}} \{e_j : j = 1, 2, \ldots \}$. Since by (2.6)
\[
S^* (E^*) \subset (\mbox{ker} S)^\bot \subset K_0 ,
\]
one obtains for every $u \in E^*$, by Lemma 2,
\[
\sum\limits_{j} \langle u, S e_j \rangle^2 = \sum\limits_{j} (S^* u, e_j)^2_K = \| S^* u \|^2_K
= \langle u, C u \rangle .
\]
The assertion follows from Lemma 1.
Conversely, if $(f_j)$ is admissible for $X$ then $(f_j)$ is a Parseval frame for $H$ by Theorem 1. Set $e_j := (S| ( \mbox{ker} S)^\bot )^{-1} (f_j) \in (\mbox{ker} S)^\perp $. Then by (2.6) and Lemma 2,
for every $k \in (\mbox{ker} S)^\bot $,
\[
\sum\limits_j (k, e_j)^2_K = \sum\limits_j (S k, f_j)2_H = || S k ||^2_H = || k ||^2_K
\]
so that again by Lemma 2, $(e_j)$ is a Parseval frame for $(\mbox{ker} S)^\bot$. \hfill{$\Box$} \\ \\
{\sc Examples} $\bullet$ Let $S : H \rightarrow E$ be the inclusion map.
Then $C = S S^*$. \\
$\bullet$ Let $K$ be the closure of $E^*$ in $L^2(\P_X)$ and $S : K \rightarrow E, Sk = \E k(X)X$. Then $S^* : E^* \rightarrow K$ is the natural embedding. Thus $C = SS^*$
and $S$ is injective (see (2.6)).
($K$ is sometimes called the energy space of $X$.) One obtains
\[
H = S(K) = \{ \E k(X)X : k \in K \}
\]
and
\[
(\E k_1 (X)X, \E k_2 (X) X)_H = \int k_1 k_2 d\P_X .
\]
$\bullet$ Let $E$ be a Hilbert space, $K = E$ and $S = C^{1/2}$. Then $C = SS^* = S2$ and
$(\mbox{ker} S)^\bot = \overline{H}$. Consequently, if $(e_j)$ is an orthonormal basis of
the Hilbert subspace $\overline{H}$ of $E$ consisting of eigenvectors of $C$ and $( \lambda_j )$ the corresponding nonzero eigenvalues, then
$( \sqrt{\lambda_j} e_j )$ is admissible for $X$ and an orthonormal basis of $(H, (\cdot , \cdot)_H)$ (Karhunen-Lo\`eve basis). \\

Admissible sequences for $X$ can be charaterized as the sequences $(S e_j)_{j \geq 1}$ where $(e_j)$ is a fixed orthonormal basis of $K$ and $S$ provides a factorization of $C$. That every sequence $(S e_j)$ of this type is admissible follows from Proposition 1.
\begin{Thm} Assume that $(f_j)_{j \geq 1}$ is admissible for $X$. Let $K$ be an infinite dimensional real separable Hilbert space and $(e_j)_{j \geq 1}$ an orthonormal basis of $K$. Then there is a factorization
$C = SS^*, S : K \rightarrow E$ such that $S(e_j) = f_j$ for every $j $.
\end{Thm}

\noindent {\bf Proof.} First, observe that $\sum\limits^\infty_{j=1} c_j f_j$ converges in $E$ for every
$(c_j)_j \in l_2 (\N)$. In fact, using Lemma 1,
\[
\begin{array}{lcl}
\| \sum\limits^{n+m}_{j=n} c_j f_j \|^2 & = & \sup_{\| u \|  \leq 1}
\langle u, \sum\limits^{n+m}_{j=n} c_j f_j \rangle^2 \\
& \leq & \sum\limits^{n+m}_{j=n} c^2_j \sup_{\| u \| \leq 1} \sum\limits^\infty_{j=1} \langle u, f_j\rangle^2\\
& = & \sum\limits^{n+m}_{j=n} c^2_j \sup_{\| u \| \leq 1} \langle u, C u \rangle \\
& = & \sum\limits^{n+m}_{j=n} c^2_j \| C \| \rightarrow 0, \; n, m \rightarrow \infty
\end{array}
\]
and thus the sequence is Cauchy in $E$. Now define $S : K \rightarrow E$ by
\[
S(k) := \sum\limits^\infty_{j=1} (k, e_j)_K f_j
\]
where $\sum (k, e_j)_K f_j$ converges in $E$ since $((k, e_j)_K)_j \in l_2 ( \N)$.
$S$ is obviously linear. Moreover, for $k \in K$, using again Lemma 1,
\[
\begin{array}{lcl}
\| S k \|^2 & = & \sup_{\| u \| \leq 1} \langle u, S k \rangle^2 \\
& = & \sup_{\|  u \| \leq 1} ( \sum\limits^\infty_{j=1} (k, e_j)_K \langle u, f_j\rangle )2 \\
& \leq & \| k \|^2_K \|  C \|  .
\end{array}
\]
Consequently,  $S$ is continuous and $S(e_j) = f_j$ for every $j$.
(At this place one needs orthonormality of $(e_j)$.)
Finally, $S^*(u) = \sum\limits^\infty_{j=1} \langle u, f_j \rangle e_j$ and hence
\[
SS^* u = \sum\limits^\infty_{j=1} \langle u, f_j \rangle f_j = C u
\]
for every $u \in E^{*}$ by Lemma 1. $\Box$ \\

It is an immediate consequence of the preceding theorem that an admissible sequence $(f_j)$ for $X$ is an orthonormal basis for $H$ if and only if $(f_j)$ is $l_2$-independent, that is
$\sum\limits^\infty_{j=1} c_j f_j = 0$ for some $(c_j) \in l_2(\N)$ implies $c_j = 0$ for every $j$. In fact, $l_2$-independence of $(f_j)$ implies that the operator $S$ in Theorem 2 is injective.

Let $F$ be a further separable Banach space and $V : E \rightarrow F$  a
$\P_X$-{\em measurable linear transfromation}, that is, $V$ is Borel measurable and linear on a Borel measurable subspace $D_V$ of $E$ with $\P_X(D_V) = 1$.
Then $H_X \subset D_V$, the operator $VJ_X : H_X \rightarrow F$ is linear and continuous, where
$J_X : H_X \rightarrow E$ denotes the inclusion map and
$V(X)$ is centered Gaussian with covariance operator
\begin{equation}
C_{V(X)} = V J_X (V J_X)^*
\end{equation}
(see [11], [2], Chapter 3.7). Consequently, by (2.6)
\begin{eqnarray}
H_{V(X)} & = & V(H_X) , \\
*[.4em] (V h_1, V h_2)_{H_{V(X)}} & = & (h_1, h_2)_{H_X} , h_1 \in H_X, h_2 \in ( \mbox{ker} (V \mid H_X))^\bot  . \nonumber
\end{eqnarray}
Note that the space of $\P_X$-measurable linear transfromation $E \rightarrow F$ is equal to the $L^p_F(\P_X)$-closure of the space of linear continuous operators $E \rightarrow F, p \in [1, \infty)$ (see [11]).

From Theorem 1 and Proposition 1 one may deduce the following proposition.
\begin{Pro}
Assume that $V : E \rightarrow F$ is a $\P_X$-measurable linear transformation. If $(f_j)_{j \geq 1}$ is admissible for $X$ in $E$, then $(V(f_j))_{j \geq 1}$ is admissible for
$V(X)$ in $F$. Conversely, if $V|H_X$ is injective and
$(g_j)_{j \geq 1}$ an admissible sequence for $V(X)$ in $F$, then there exists a sequence $(f_j)_{j \geq 1}$ in $E$ which is admissible for $X$ such that $V(f_j) = g_j$ for every $j$.
\end{Pro}
{\sc Example}
Let $X$ and $Y$ be jointly centered Gaussian random vectors in $E$ and $F$, respectively. Then $\E(Y|X) = V(X)$ for some $\P_X$-measurable linear transformation $V : E \rightarrow F$.
The cross covariance operator $C_{YX} : E^{*} \rightarrow F$,
$C_{YX} u = \E \langle u, x \rangle Y$ can be factorized as
$C_{YX} = U_{YX} S^{*}_X$, where $C_X = S_X S^{*}_X$ is the energy factorization of $C_X$ with $K_X$ the closure of $E^{*}$ in $L^2(\P_X)$ and $U_{YX} : K_X \rightarrow F, U_{YX} k = \E k(X)Y$.
Then
\[
V = U_{YX} S^{-1}_X \; \mbox{on} \; H_X
\]
(see [11]). Consequently, if $(f_j)_{j \geq 1}$ is admissible for $X$ in $E$ then
$(U_{YX} S^{-1}_X f_j)_{j \geq 1}$ is admissible for $\E(Y|X)$ in $F$.
\section{Continuous Gaussian processes}
\setcounter{equation}{0}
\setcounter{Assumption}{0}
\setcounter{Theorem}{0}
\setcounter{Proposition}{0}
\setcounter{Corollary}{0}
\setcounter{Lemma}{0}
\setcounter{Definition}{0}
\setcounter{Remark}{0}
Now let $I$ be a compact metric space and $X = (X_t)_{t \in I}$ be a real pathwise continuous centered Gaussian process. Let $E := {\cal C}(I)$ be equipped with the sup-norm
$\|x\| =\sup_{t \in I} |x(t)|$ so that the norm dual ${\cal C}(I)^{*}$ coincides with the space of finite signed Borel measures on $I$ by the Riesz theorem. Then $X$ can be seen as a ${\cal C}(I)$-valued Gaussian random vector and the covariance operator $C : {\cal C}(I)^* \rightarrow {\cal C}(I)$ takes the form
\begin{eqnarray}
C u(t) & = & \langle \delta_t , C u \rangle = \langle C \delta_t , u \rangle \nonumber \\
 & = & \langle \E X_t X, u \rangle = \int_I \E X_t X_s d u (s) .
\end{eqnarray}
\begin{Cor}
Let $(f_j)_{j \geq 1} \in {\cal C} (I)^\N$. \\
(a) If
\[
\E X_s X_t = \sum\limits^\infty_{j=1} f_j (s) f_j(t) \; \mbox{for every} \; s, t \in I
\]
then $(f_j)$ is admissible for $X$. \\ \\
(b) If
\[
\sum\limits^\infty_{j=1} f_j(t)2 < \infty \; \mbox{for every} \; t \in I
\]
and if the process $Y$ with $Y_t = \sum\limits^\infty_{j=1} \xi_j f_j(t)$ has a pathwise continuous modification $X$, then $(f_j)$  is admissible for $X$ and $X = \sum\limits^\infty_{j=1} \xi_j f_j$ a.s.
\end{Cor}
{\bf Proof.}
(a) For $u \in G := \mbox{span} \; \{ \delta_t : t \in I\}, u = \sum\limits^m_{i=1} \alpha_i \delta_{t_i}$ we have
\[
\langle u, C u \rangle = \sum\limits^m_{i=1} \sum\limits^m_{k=1} \alpha_i \alpha_k \E X_{t_i} X_{t_k}
\]
and
\[
\sum\limits^n_{j=1} \langle u, f_j \rangle^2 = \sum\limits^m_{i=1} \sum\limits^m_{k=1} \alpha_i \alpha_k
\sum\limits^n_{j=1} f_j(t_i) f_j(t_k)
\]
so that
\[
\sum\limits^\infty_{j=1} \langle  u, f_j \rangle^2 = \langle u, C u\rangle .
\]
Since $G$ is a separating subspace of ${\cal C}(I)^{*}$ the assertion follows from Lemma 1. \\
(b) Notice that $\sum \xi_j f_j(t)$ converges a.s. in $\R$ and $Y$ is a centered Gaussian process. Hence $X$ is centered Gaussian. Since
\[
\E X_s X_t = \E Y_s Y_t = \sum\limits^\infty_{j=1} f_j (s) f_j(t)\; \mbox{for every} \; s, t \in I,
\]
the assertion follows from (a). \hfill{$\Box$} \\

Factorizations of $C$ can be obtained as follows.
For Hilbert spaces $K_i$, let $\oplus^m_{i=1} K_i$ denote the Hilbertian (or $l_2-$)direct sum.
\begin{Lem} For $i \in \{ 1 , \ldots, m\}$, let $K_i$ be a real separable Hilbert space. Assume the representation
\[
\E X_s X_t = \sum\limits^m_{i=1} (g^i_s, g^i_t)_{K_i} , s, t \in I
\]
for vectors $ g^i_ t \in K_i$. Then
\[
S : \oplus^m_{i=1} K_i \rightarrow {\cal C}(I), S k (t) := \sum\limits^m_{i=1} (g^i_t, k_i)_{K_i}
\]
is a linear continuous operator,
$(\mbox{ker} S)^\bot = \overline{span} \{(g1_{t}, \ldots , g^m_t) : t \in I\}$ and $C = SS^*$.
\end{Lem}
{\bf Proof.}
Let $K := \oplus^m_{i=1} K_i$ and $g_t := (g1_{t}, \ldots , g^m_t)$. Then
$\E X_s X_t = (g_s , g_t)_K$ and $Sk(t) = (g_t, k)_K$.
First, observe that
\[
\sup_{t \in I} \| g_t \|_K \leq \| C \|^{1/2} < \infty .
\]
Indeed, for every $t \in I$, by (3.1),
\[
\| g_t \|^2_K = \E X2_t = \langle \delta_t, C \delta_t \rangle \leq \| C \|  .
\]
The function $S k$ is continuous for $k \in \; \mbox{span} \; \{ g_s : s \in I \}$. This easily implies that $S k$ is continuous for every
$k \in \; \overline{\mbox{span}} \{ g_s : s \in I \}$ and thus for every $k \in K$. $S$ is obviously linear and
\[
\| S k\|  = \sup_{t \in I} \mid (g_t, k)_K \mid \leq \| C \|^{1/2} \| k \|_K .
\]
Finally, $S^*(\delta_t) = g_t$ so that
\[
SS^* \delta_t (s) = S g_t (s) = \E X_s X_t = C \delta_t (s)
\]
for every $s, t \in I$. Consequently, for every $u \in {\cal C}(I)^*, t \in I$,
\[
\begin{array}{lcl}
SS^* u(t) & = & \langle SS^* u , \delta_t\rangle = \langle u, SS^* \delta_t \rangle \\
& = & \langle u, C \delta_t \rangle = \langle C u, \delta_t \rangle = C u (t)
\end{array}
\]
and hence $C = SS^*$. $\Box$%

\bigskip
\noindent {\sc Example} Let $K$ be the first Wiener chaos, that is
$K = \; \overline{\mbox{span}}  \{ X_t : t \in I \}$ in $L^2 (\P)$ and $g_t = X_t$.
Then $S k = \E \,k X$ and $S$ is injective. If for instance $X = W$ (Brownian motion) and $I = [0,T]$, then
\[
K = \left\{ \int^T_0 f (s) dW_s : f \in L^2 ([0,T], dt) \right\}.
\]

We derive from the preceeding lemma and Proposition 1 the following corollary.
\begin{Cor} Assume the situation of Lemma 3. Let $(e^i_j)_j $ be a Parseval frame sequence in $K_i$ satisfying
$\{g^i_t : t \in I\} \subset \; \overline{span} \; \{ e^i_j : j = 1, 2, \ldots \}$. Then,
$(S_i(e^i_j))_{1 \leq i \leq m, j}$
is admissible for $X$, where $S_i k(t) = (g^i_t, k)_{K_i}$.
\end{Cor}

The next corollary implies the well known fact that the Karhunen-Lo\`eve expansion of $X$ in some Hilbert space $L^2 (I, \mu)$ already converges uniformly in $t \in I$. It appears as special case of Proposition 2.
\begin{Cor} Let $\mu$ be a finite Borel measure on $I$ with $\mbox{supp} (\mu) = I$ and let $V : {\cal C}(I) \rightarrow L^2 (I, \mu)$ denote the natural (injective) embedding. Let $(g_j)_{j \geq 1}$ be admissible for $V(X)$ in $L^2(I, \mu)$. Then there exists a sequence
$(f_j)_{j \geq 1}$ in ${\cal C}(I)$ which is admissible for $X$ such that $V(f_j) = g_j$ for every $j$.
\end{Cor}

The admissibility feature is stable under tensor products. For $i \in \{1 , \ldots , d \}$, let $I_i$
be a compact metric space and $X^i = (X^i_t)_{t \in I_i}$ a continuous centered Gaussian process.
Set $I := \Pi^d_{i=1} I_i$ and let $X = (X_t)_{t \in I}$ be a continuous centered Gaussian process with
covariance function
\begin{equation}
E X_s X_t = \Pi^d_{i=1} \E X^i_{s_i} X^i_{t_i} , s, t \in I .
\end{equation}
For instance, $X := \otimes^d_{i=1} X^i$ satisfies (3.2) provided $X1, \ldots , X^d$ are independent. For real separable Hilbert spaces $K_i$, let $\widehat{\otimes^d_{i=1}} K_i$ denote the $d$-fold Hilbertian tensor product.
\begin{Pro} For $i \in \{1, \ldots , d \}$, let $(f^i_j)_{j \geq 1}$ be an admissible sequence for
$X^i$ in ${\cal C}(I_i)$. Then
\[
( \otimes ^d_{i=1} f^i_{j_i} )_{ \underline{j} = (j_1 , \ldots , j_d) \in \N^d }
\]
is admissible for $X$ with covariance (3.2) in ${\cal C}(I)$. Furthermore, if $C_{X^i} = S_i S^{*}_i, S_i : K_i \rightarrow {\cal C} (I_i)$
is a factorization of $C_{X^i}$, then
$\otimes^d_{i=1} S_i : \widehat{\otimes ^d_{i=1} } K_i \rightarrow {\cal C}(I)$ provides a factorization of $C_X$.
\end{Pro}
{\bf Proof.} For $i \in \{1, \ldots , d\}$, let $K_i$ be a real separable Hilbert space and
$(e^i_j)_{j}$ an orthonormal basis of $K_i$. Then $( \otimes^d_{i=1} e^i_{j_i} )_{\underline{j}}$ is an orthonormal basis of
$K := \widehat{\otimes^{d}_{i=1}} K_i$.

If $C_{X^i} = S_i S^{*}_i, S_i : K_i \rightarrow {\cal C} (I_i)$ is a factorization of $C_{X^i}$, set $g^i_t := S^{*}_i \delta_t, t \in I_i$. Then
$\E X^i_s X^i_t = ( g^i_s, g^i_t)_{K_i}$ and hence, by (3.2)
\[
\E X_s X_t = \Pi^d_{i=1} (g^i_{s_i}, g^i_{t_i})_{K_i} = ( \otimes^d_{i=1} g^i_{s_i} ,
\otimes^d_{i=1} g^i_{t_i} )_K , s, t \in I.
\]
Consequently, by Lemma 3
\[
U : K \rightarrow {\cal C}(I), U k(t) = ( \otimes^d_{i=1} g^i_{t_i} , k)_K
\]
provides a factorization of $C_X$. Since
\[
\begin{array}{lcl}
U( \otimes^d_{i=1} e^i_{j_i} )(t) & = & \Pi^d_{i=1} (g^i_{t_i}, e^i_{j_i})_{K_i} \\
& = & \Pi^d_{i=1} S_i e^i_{j_i} (t_i) = \otimes ^d_{i=1} (S_i e^i_{j_i} )(t) \\
& = & ( \otimes^d_{i=1} S_i) ( \otimes^d_{i=1} e^i_{j_i} )(t) , t \in I,
\end{array}
\]
we obtain $U = \otimes^d_{i=1} S_i$ and thus $\otimes^d_{i=1} S_i$ provides a factorization of
$C_X$.

If $(f^i_j)_{j \geq 1}$ is admissible for $X^i$, then by Theorem 2 assuming now that $K_i$ is infinite dimensional, there is a factorization
$C_{X^i} = T_i T^{*}_i , T_i : K_i \rightarrow {\cal C}(I_i)$ such that $T_i(e^i_j) = f^i_j$ for every $j$.
Since $\otimes^d_{i=1} T_i : K \rightarrow {\cal C}(I)$ provides a factorization of $C_X$ as shown above and
$(\otimes^d_{i=1} T_i )(\otimes^d_{i=1} e^i_{j_i}) = \otimes^d_{i=1} f^i_{j_i}$,
it follows from Proposition 1 that $( \otimes^d_{i=1} f^i_{j_i})_{\underline{j} \in \N^d}$ is admissible for $X$. $\Box$ \\ \\
{\bf Comments on optimal expansions.} For $n \in \N$, let
\begin{equation}
l_n(X) := \inf \{ \E || \sum\limits^\infty_{j=n} \xi_j f_j || : (f_j)_{j \geq 1} \in {\cal C} (I)^\N \; \mbox{admissible for} \; X \} .
\end{equation}
Rate optimal solutions of the $l_n(X)$-problem are admissible sequences $(f_j)$ for $X$ in ${\cal C}(I)$ such that
\[
\E || \sum\limits^\infty_{j=n} \xi_j f_j || \approx l_n (X) \; \mbox{as} \; n \rightarrow \infty.
\]

For $I = [0,T]^d \subset \R^d$, consider the covariance operator $R = R_X$ of $X$ on $L^2(I,dt)$ given by
\begin{equation}
R : L^2(I,dt) \rightarrow L^2(I,dt), Rk(t) = \int_I \E X_sX_t k(s) ds.
\end{equation}
Using (3.1) we have $R_X = V C_X V^{*}$, where $V : C(I) \rightarrow L^2(I,dt)$
denotes the natural (injective) embedding. The choice of Lebesgue measure on I is the best choice for our purposes (see (A1)). Let $\lambda_1 \geq \lambda_2 \geq \ldots > 0$ be the ordered nonzero eigenvalues of $R$ (each written as many times as its multiplicity).
\begin{Pro}
Let $I = [0,T]^d$. Assume that the eigenvalues of $R$ satisfy \\ \\
(A1)
$\lambda_j \geq c_1 j^{-2 \vartheta} \log (1+j)^{2\gamma}$ for every $j \geq 1$ with $\vartheta > 1/2, \gamma \geq 0$ and $c_1 > 0$ \\ \\
and that $X$ admits an admissible sequence $(f_j)$ in ${\cal C}(I)$ satisfying \\ \\
(A2)
$|| f_j|| \leq c_2 j^{-\vartheta} \log (1+j)^\gamma$ for every $j \geq 1$with $c_2 < \infty$, \\ \\
(A3)
$f_j$ is a-H{\"o}lder-continuous and $[f_j]_a \leq c_3 j^b$ for every $j \geq 1$ with
$a \in (0,1], b \in \R$ and $c_3 < \infty$, where
\[
[f]_a = \sup_{s \not= t} \frac{|f(s) - f(t)|}{|s-t|^a}
\]
(and $|t|$ denotes the $l_2$-norm of $t \in \R^d$). \\
Then
\begin{equation}
l_n(X) \approx n^{-(\vartheta - \frac{1}{2} ) } (\log n)^{\gamma + \frac{1}{2} } \; \mbox{as} \; n \rightarrow \infty
\end{equation}
and $(f_j)$ is rate optimal.
\end{Pro}
{\bf Proof.}
The lower estimate in (3.5) follows from (A1) (see [8], Proposition 4.1) and from (A2) and (A3) follows
\[
\E || \sum\limits^\infty_{j=n} \xi_j f_j || \leq c_4 n^{-(\vartheta - \frac{1}{2})} ( \log
( 1+ n))^{\gamma + \frac{1}{2} }
\]
for every $n \geq 1$, (see [13], Theorem 1). \hfill{$\Box$} \\

Concerning assumption (A3) observe that we have by (2.2) and (3.1) for
$h \in H = H_X, s, t \in I, $
\[
h(t) = < \delta_t, h  >  = ( h, C \delta_t)_H
\]
and
\[
|| C(\delta_s - \delta_t)||^2_H = < \delta_s - \delta_t , C ( \delta_s - \delta_t)> = \E | X_s - X_t|^2
\]
so that
\begin{eqnarray}
|h(s) - h(t)| & = & |(h, C(\delta_s - \delta_t))_H| \nonumber \\
& \leq & || h ||_H || C ( \delta_s - \delta_t)||_H \nonumber \\
& = & || h ||_H ( \E |X_s - X_t|^2)^{1/2} .
\end{eqnarray}

Consequently, since admissible sequences are contained in the unit ball of $H$, (A3) is satisfied with $b = 0$ provided $I \rightarrow L^2(\P), t \mapsto X_t$ is $a$-H{\"o}lder-continuous.

The situation is particularly simple for Gaussian sheets.
\begin{Cor}
Assume that for $i \in \{ 1, \ldots , d\}$, the continuous centered Gaussian process
$X^i = (X^i_t)_{t \in [0,T]}$ satisfies (A1) - (A3) for some admissible sequence
$(f^i_j)_{j \geq 1}$ in ${\cal C}([0,T])$ with parameters $\vartheta_i, \gamma_i, a_i, b_i$ such that $\gamma_i = 0$ and let $X = (X_t)_{t \in I}, I = [0,T]^d$ be the continuous centered Gaussian sheet with covariance (3.2). Then
\begin{equation}
l_n(X) \approx n^{-(\vartheta - \frac{1}{2})} ( \log n)^{\vartheta(m-1) + \frac{1}{2} }
\end{equation}
with
$\vartheta = \min_{1 \leq i \leq d} \vartheta_i$ and $m = \mbox{card} \{ i \in \{ 1 , \ldots , d \}: \vartheta_i = \vartheta \}$
and a decreasing arrangement of $(\otimes^d_{i=1} f^i_{j_i})_{{\underline j} \in \N^d}$
is rate optimal for $X$.
\end{Cor}
{\bf Proof.} In view of Lemma 1 in [13] and Proposition 3, the assertions follow from Proposition 4. \hfill{$\Box$} \\ \\
{\sc Examples} The subsequent examples may serve as illustrations. \\
$\bullet$ Let
$W =(W_t)_{t \in [0,T]}$ be a {\em standard Brownian motion}. Since
$\E W_s W_t = s \wedge t = \int^T_0 1_{[0,s]} (u) 1_{[0,t]}(u) du$, the (injective) operator
\[
S : L^2 ([0,T], dt) \rightarrow {\cal C} ([0,T]), \quad Sk(t) = \int^t_0 k(s) ds
\]
provides a factorization of $C_W$ so that we can apply Corollary 2.
The orthonormal basis $e_j (t) = \sqrt{2/T} \cos ( \pi (j - 1/2)t/T), j \geq 1$ of
$L^2 ([0,T], dt)$ yields the admissible sequence
\begin{equation}
f_j (t) = S e_j (t) = \frac{ \sqrt{2 T}}{\pi(j-1/2) } \sin ( \frac{\pi (j-1/2)t}{T} ) ,\; j \geq 1
\end{equation}
for $W$ (Karhunen-Lo\`eve basis of $H_W$) and $e_j (t) = \sqrt{2/T} \sin ( \pi j t /T)$ yields the admissible sequence
\[
g_j(t)  =  \frac{\sqrt{2T}}{\pi j} (1 - \cos ( \frac{\pi j t}{T})),\; j \geq 1 .
\]
Then
\begin{eqnarray}
f1_j (t) = \frac{1}{\sqrt{2}} f_j(t) & = & \frac{\sqrt{T}}{\pi(j-1/2)} \sin ( \frac{\pi (j-1/2)t}{T} ) , \;j \geq 1 \\
f2_j(t) = \frac{1}{\sqrt{2}} g_j(t)  & = & \frac{\sqrt{T}}{\pi j} (1 - \cos ( \frac{\pi j t}{T})), \;j \geq 1 \nonumber
\end{eqnarray}
is a Parseval frame for $H_W$ and hence admissible for $W$.
The trigonometric basis
$e_0 (t) = 1/\sqrt{T}$,
$e_{2j} (t) = \sqrt{2/T} \cos(2 \pi j t/T), e_{2j-1} (t) = \sqrt{2/T} \sin (2 \pi j t /T)$
of $L^2([0,T], dt)$ yields the admissible sequence
\begin{eqnarray}
f_0 (t) = \frac{t}{\sqrt{T}} , f_{2 j} (t) = \frac{\sqrt{T}}{\sqrt{2} \pi j} \sin (
\frac{2 \pi j t}{T}) , \\
f_{2 j-1} (t) = \frac{\sqrt{T}}{ \sqrt{2} \pi j} (1 - \cos ( \frac{2 \pi j t}{T})), j \geq 1 \nonumber
\end{eqnarray}
(Paley-Wiener basis of $H_W$). By Proposition 4, all these admissible sequences for $W$ (with $f_{2j} := f1_j, f_{2j-1} := f^2_j$, say in (3.9)) are rate optimal.

Assume that the wavelet system $2^{j/2} \psi(2^j \cdot - k), j, k \in \Z$
is an orthonormal basis (or only a Parseval frame) for $L^2(\R, dt)$. Then the restrictions of these functions to [0,T] clearly provide a Parseval frame for
$L^2([0,T], dt)$ so that the sequence
\[
f_{j,k}(t) = S(2^{j/2} \psi(2^j \cdot -k))(t) = 2^{-j/2} \int^{2^jt-k}_{-k}
\psi(u) du, j, k \in \Z
\]
is admissible for $W$. If $\psi \in L1(\R, dt)$ and
$\Psi(x) := \int^x_{- \infty} \psi(u) du$, then this admissible sequence takes the form
\begin{equation}
f_{j,k} (t) = 2^{-j/2} (\Psi(2^j t-k) - \Psi(-k)), j, k \in \Z .
\end{equation}

\ms
\noindent $\bullet$ We consider the Dzaparidze-van Zanten expansion of the {\em fractional Brownian motion}
$X = (X_t)_{t \in [0,T]}$ with Hurst index $\rho \in (0,1)$ and covariance function
\[
\E X_s X_t = \frac{1}{2} (s^{2\rho} + t^{2\rho} - \mid s-t \mid^{2\rho} ) .
\]
These authors discovered in [5] for $T = 1$ a time domain representation
\[
\E X_s X_t = (g^1_s, g^1 _t)_K + (g^2_s , g2_t)_K
\]
with $K = L^2 ([0,1], dt)$ and kernels $g^i_t \in L^2 ([0,1], dt)$. Hence by Lemma 3, the operator
\[
S : L^2([0,1], dt) \oplus L^2 ([0,1], dt) \rightarrow {\cal C} ([0,1]), \;S(k_1, k_2)(t) =
\int^1_0 g^1_t (s) k_1 (s) ds + \int1_0 g^2_t (s) k_2 (s) ds
\]
provides a factorization of $C_X$ so that for every pair of orthonormal bases
$(e1_j)_{j \geq 1}$ and $(e^2_j)_{j \geq 1}$ of $L^2 ([0,1], dt)$,
\[
f^i_j (t) = \int^1_0 g^i_t (s) e^i_j (s) ds , \;j \geq 1,\; i = 1,2
\]
is admissible in ${\cal C}([0,1])$ for $X$. By Corollary 2, this is a consequence of the above representation of the covariance function~(and needs no extra work). Then Dzaparidze and van Zanten [5] could calculate $f^i_j$ explicitely for the Fourier-Bessel basis of order $- \rho$ and $1-\rho$, respectively and arrived at the admissible family in ${\cal C}([0,1])$
\[
\begin{array}{lcl}
f^1_j (t) & = &
\displaystyle{ \frac{ c_\rho \sqrt{2} }{ \mid J_{1 -\rho} (x_j) \mid x^{\rho+1}_j } }
\sin (x_j t) , j \geq 1 \\
f2_j (t) & = &
\displaystyle{\frac{c_\rho \sqrt{2} }{ \mid J_{-\rho} (y_j) \mid y^{\rho+1}_j }  }
(1 - \cos (y_j t)) , j \geq 1 \nonumber
\end{array}
\]
where $J_\nu $ denotes the Bessel function of the first kind of order $\nu$, $0 < x_1 < x_2 < \ldots$ are the positive zeros of $J_{-\rho} , 0 < y_1 < y_2 < \ldots $ the positive zeros of $J_{1 - \rho}$ and
$c^2_\rho = \Gamma (1 + 2 \rho) \sin (\pi \rho)/\pi$.
Consequently, by self-similarity of $X$, the sequence
\begin{eqnarray}
f^1_j (t) & = & \frac{T^\rho c_\rho \sqrt{2}}{ \mid J_{1-\rho} (x_j) \mid  x^{\rho+1}_j } \sin ( \frac{x_jt}{T} ), j \geq 1 \\
f2_j (t) & = & \frac{T^\rho c_\rho \sqrt{2}}{ \mid J_{-\rho} (y_j) \mid  y^{\rho+1}_j } (1 - \cos ( \frac{y_jt}{T})), j \geq 1 \nonumber
\end{eqnarray}
in ${\cal C}([0,T])$ is admissible for $X$.
Using Lemma 1, one can deduce (also without extra work)
\[
\E X_s X_t = \sum\limits^\infty_{j=1} f^1_j (s) f^1_j(t) + \sum\limits^\infty_{j =1 } f^2_j (s) f2_j(t)
\]
uniformly in $(s,t) \in [0,T]^2$. Rate optimality of (3.12) (using an arrangement like
$f_{2j} := f^1_j$,  $f_{2j-1} := f2_j)$ is shown in [6] based on the work [8]
and is also an immediate consequence of Proposition 4 since
\[
x_j \sim y_j \sim \pi j, J_{1-\rho} (x_j) \sim J_{-\rho} (y_j) \sim \frac{\sqrt{2}}{\pi}
j^{-1/2}
\]
(see [5]), and the eigenvalues satisfy $\lambda_j \sim c j^{-(1 + 2 \rho)}$ as $j \rightarrow \infty$ (see [3], [12]).

In the ordinary Brownian motion case $\rho = 1/2,$ (3.12) coincides with (3.9).
The interesting extension of (3.10) to fractional Brownian motions is discussed in [7]
and extensions of the wavelet expansion (3.11) can be found in [1], [14]. \\
$\bullet$ Let $X = (X_t)_{t \in [0,T]}$ be {\em Brownian bridge} with covariance
\[
\E X_s X_t = s \wedge t - \frac{st}{T} = \int^T_0 (1_{[ 0,s]} (u) - \frac{s}{T})(1_{[0,t]}(u) - \frac{t}{T}) du .
\]
By Lemma 3, the operator
\[
S : L^2 ([0,T], dt) \rightarrow {\cal C} ([0,T]), Sk(t) = \int^t_0 k (s) ds - \frac{t}{T} \int^T_0 k(s) ds
\]
provides a factorization of $C_X$ and ker$S = \mbox{span} \{ 1_{[0,T]}\}$. The choice $e_j (t) = \sqrt{2/T} \cos (\pi j t/T), j \geq 1$ of an orthonormal basis of (ker $S)^\perp $  yields admissibility of
\begin{equation}
f_j (t) = S e_j (t) = \frac{\sqrt{2 T}}{\pi j} \sin ( \frac{\pi j t}{T}), j \geq 1
\end{equation}
for $X$ (Karhunen-Lo\`eve basis of $H_X$). By Proposition 4, this sequence is rate optimal. \\
$\bullet$ One considers the {\em stationary Ornstein-Uhlenbeck process} as the solution of the Langevin equation
\[
d X_t = - \alpha  X_t dt + \sigma d W_t , t \in [0,T]
\]
with $X_0$ independent of $W$and $N(0, \frac{\sigma2}{2 \alpha}$)-distributed, $\sigma > 0, \alpha > 0$. It admits the explicit representation
\[
X_t = e^{- \alpha t} X_0 + \sigma e^{-\alpha t} \int^t_0 e^{\alpha s} d W_s
\]
and
\[
\E X_s X_t = \frac{\sigma^2}{2 \alpha} e^{-\alpha \mid  s-t \mid } = \frac{\sigma^2}{2 \alpha} e^{-\alpha(s+t)} + \sigma2 e^{-\alpha(s+t)} \int^{s \wedge t}_0 e^{2 \alpha u} du .
\]
Thus  the (injective) operator
\[
S : \R \oplus L^2 ([0,T], dt) \rightarrow {\cal C} ([0,T]), S (c,k)(t) = \frac{c \sigma}{\sqrt{2 \alpha}} e^{-\alpha t} +
\sigma \int^t_0 e^{- \alpha(t-s)} k(s) ds
\]
provides a factorization of $C_X$ so that for every Parseval frame
$(e_j)_{j \geq 1}$ for $L^2([0,T], dt)$, the functions
\begin{equation}
f_0 (t) = \frac{\sigma}{\sqrt{2 \alpha}} e^{-\alpha t} ,\; f_j (t) = \sigma \int^t_0 e^{-\alpha(t-s)} e_j (s) ds, j \geq 1
\end{equation}
provide an admissible sequence for $X$. For instance the choice of the orthonormal basis $\sqrt{2/T} \cos (\pi(j - 1/2) t/T), j \geq 1$ implies that (3.14) is rate optimal. This follows from Lemma 1 in [13] and Proposition 4.

Another representation is given by the Lamperti transformation $X = V(W)$ for the linear continuous operator
\[
V : {\cal C} ([0, e^{2 \alpha T}]) \rightarrow {\cal C} ([0,T], V x (t) = \frac{\sigma}{\sqrt{2\alpha}} e^{-\alpha t} x(e^{2 \alpha t}) .
\]
The admissible sequence $(f_j)$ in ${\cal C} ([0, e^{2 \alpha T} ])$ for $(W_t)_{t \in [0, e^{2 \alpha T}]}$
from (4.1) yields the admissible sequence
\begin{equation}
\tilde{f}_j (t) = Vf_j (t) = \frac{\sigma}{\sqrt{\alpha} \pi (j - 1/2)} e^{\alpha(T-t)} \sin ( \pi (j-1/2) e^{-2\alpha(T-t)} ), j \geq 1
\end{equation}
for $X$. By Proposition 4, the sequence (3.15) is rate optimal.\\
$\bullet$ Sheet versions can be deduced from Proposition 3 and Corollary 4 (and need no extra work).
\section{Optimal expansion of fractional Ornstein-Uhlenbeck processes}
\setcounter{equation}{0}
\setcounter{Assumption}{0}
\setcounter{Theorem}{0}
\setcounter{Proposition}{0}
\setcounter{Corollary}{0}
\setcounter{Lemma}{0}
\setcounter{Definition}{0}
\setcounter{Remark}{0}
The fractional Ornstein-Uhlenbeck process $X^\rho = ( X^\rho_t)_{t \in \R}$ of index $\rho \in (0,2)$ is a continuous stationary centered Gaussian process having the covariance function
\begin{equation}
\E X^\rho_s X^\rho_t = e^{- \alpha | s - t |^\rho} , \alpha > 0 .
\end{equation}
We derive explicit optimal expansions of $X^\rho$ for $\rho \leq 1$. Let
\[
\gamma^\rho : \R \rightarrow \R, \gamma^\rho (t) = e^{-\alpha|t|^\rho}
\]
and for a given $T > 0$, set
\begin{equation}
\beta_0 (\rho) := \frac{1}{2 T} \int^T_{-T} \gamma^\rho(t) dt , \beta_j ( \rho) := \frac{1}{T} \int^T_{-T}
\gamma^\rho(t) \cos ( \pi j t/T) dt, j \geq 1 .
\end{equation}
\begin{Thm}
Let $\rho \in (0,1]$. Then $\beta_j ( \rho) > 0$ for every $j \geq 0$ and the sequence
\begin{eqnarray}
f_0 & = & \sqrt{\beta_0 ( \rho) }, f_{2j} = \sqrt{\beta_j ( \rho)} \cos ( \pi j t/T) , \\
f_{2 j - 1} (t) & = & \sqrt{\beta_j (\rho) } \sin ( \pi j t/T), j \geq 1 \nonumber
\end{eqnarray}
is admissible for $X^\rho$ in ${\cal C} ([0,T])$. Furthermore,
\[
l_n ( X^\rho) \approx n^{-\rho/2} ( \log n)^{1/2} \; \mbox{as} \; n \rightarrow \infty
\]
and the sequence (4.3) is rate optimal.
\end{Thm}
{\bf Proof.}
Since $\gamma^\rho$ is of bounded variation (and continuous) on [-T,T], it follows from the Dirichlet criterion that its (classical) Fourier series converges pointwise to $\gamma^\rho$ on [-T,T], that is using symmetry of $\gamma^\rho$,
\[
\gamma^\rho(t) = \beta_0(\rho) + \sum\limits^\infty_{j=1} \beta_j ( \rho) \cos ( \pi j t/T), t \in [-T,T].
\]
Thus one obtains the representation
\begin{equation}
\E X_s X_t = \gamma^\rho (s-t) = \beta_0 (\rho) + \sum\limits^\infty_{j=1} \beta_j ( \rho) [ \cos ( \pi j s /T) \cos ( \pi j t/T) + \sin ( \pi j s / T) \sin ( \pi j t/T)], s, t \in [0,T] .
\end{equation}
This is true for every $\rho \in (0,2)$. If $\rho = 1$, then integration by parts yields
\begin{equation}
\beta_0 (1) = \frac{1 - e^{-\alpha T}}{\alpha T} , \beta_j (1) =
\frac{2 \alpha T (1 - e^{-\alpha T} (-1)^j) }{ \alpha^2 T^2 + \pi^2 j2} , j \geq 1 .
\end{equation}
In particular, we obtain $\beta_j(1) > 0$ for every $j \geq 0$. If
$\rho \in (0,1)$, then $\gamma^\rho | [0, \infty)$ is the Laplace transform of a suitable one-sided strictly $\rho$-stable distribution with Lebesgue-density $q_\rho$. Consequently, for $j \geq 1$
\begin{eqnarray}
\beta_j (\rho) & = & \frac{2}{T} \int^T_0 e^{-\alpha t^\rho} \cos ( \pi j t/T) dt \\
    & = & \int^\infty_0 \frac{2}{T} \int^T_0 e^{- t x} \cos ( \pi j t/T) dt q_\rho (x) dx \nonumber \\
    & = & \int^\infty_0 \frac{ 2 x T (1 - e^{-xT} (-1)^j) }{ x^2 T^2 + \pi^2 j2 }
    q_\rho (x) dx . \nonumber
\end{eqnarray}
Again, $\beta_j(\rho) > 0$ for every $j \geq 0$. It follows from (4.4) and Corollary 1(a) that the sequence $(f_j)_{j \geq 0}$ defined in (4.3) is admissible for $X^\rho$ in ${\cal C}([0,T])$.

Next we investigate the asymptotic behaviour of $\beta_j ( \rho)$ as $j \rightarrow \infty$ for
$\rho \in (0,1)$. The spectral measure of $X^\rho$ still for $\rho \in (0,2)$ is a symmetric $\rho$-stable distribution with continuous density $p_\rho$ so that
\[
\begin{array}{lcl}
 \gamma^\rho(t) & = & \displaystyle\int_\R e^{itx} p_\rho (x) dx \\
  & = & \displaystyle 2 \int^\infty_0 \cos (tx) p_\rho (x) dx , t \in \R
\end{array}
\]
and the spectral density satisfies the high-frequency condition
\begin{equation}
p_\rho (x) \sim c (\rho) x^{-(1+\rho)} \; \mbox{as} \; x \rightarrow \infty
\end{equation}
where
\[
c ( \rho) = \frac{\alpha \Gamma (1+\rho) \sin ( \pi \rho/2) }{ \pi } .
\]
Since by the  Fourier inversion formula
\[
p_\rho (x) = \frac{1}{\pi} \int^\infty_0 \gamma^\rho (t) \cos (tx) dt, x \in \R ,
\]
we obtain for $j \geq 1$,
\[
\begin{array}{lcl}
\beta_j ( \rho) & = & \frac{2}{T} \int^T_0 \gamma^\rho(t) \cos ( \pi j t / T) dt \\ \\
& = & \frac{2}{T} ( \int^\infty_0 \gamma^\rho (t) ) \cos ( \pi j t/T) dt - \int^\infty_T \gamma^\rho (t) \cos ( \pi j t /T) dt) \\ \\
& = & \frac{2 \pi}{T} p_\rho (\pi j/T) - \frac{2}{T} \int^\infty_T \gamma^\rho (t) ) \cos ( \pi j t/T) dt \end{array}
\]
Integrating twice by parts yields
\[
\int^\infty_T \gamma^\rho (t) \cos ( \pi j t/T) dt = O (j^{-2})
\]
for any $\rho \in (0,2)$ so that for $\rho \in (0,1)$
\begin{eqnarray}
\beta_j (\rho) & \sim & \frac{2 \pi}{T} p_\rho ( \pi j/T) \\
    & \sim & \frac{ 2 \pi T^\rho c ( \rho) }{ ( \pi j)^{1+ \rho} } \nonumber \\
    & = & \frac{ 2 \alpha T^\rho \Gamma (1+\rho) \sin ( \pi \rho/2) }{ (\pi j)^{1+\rho} } \;
    \mbox{as} \; j \rightarrow \infty . \nonumber
\end{eqnarray}

We deduce from (4.5) and (4.8) that the admissible sequence (4.3) satisfies the conditions (A2) and (A3) from Proposition 4 with parameters $\vartheta = (1 + \rho)/2, \gamma = 0, a = 1$ and
$b = (1-\rho)/2$. Furthermore, by Theorem 3 in Rosenblatt [16] the asymptotic behaviour of the eigenvalues of the covariance operator of
$X^\rho$ on $L^2 ([0,T], dt)$ (see (3.4)) for $\rho \in (0,2)$ is as follows:
\begin{equation}
\lambda_j \sim \frac{2 T^{1+\rho} \pi c(\rho) }{ (\pi j )^{1+\rho} }
\; \mbox{as} \; j \rightarrow \infty .
\end{equation}
Therefore, the remaining assertions follow from Proposition 4. \hfill{$\Box$} \\

Note that the admissible sequence (4.3) is not an orthonormal basis for $H = H_{X^\rho}$ but only a Parseval frame at least in case $\rho = 1$. In fact, it is well known that for $\rho = 1$,
\[
||h||^2_H = \frac{1}{2} ( h(0)^2 + h (T)^2 ) + \frac{1}{2\alpha} \int^T_0 ( h^{'} (t)^2 + \alpha2 h (t)^2 ) dt
\]
so that e.g.
\[
|| f_{2 j - 1} ||^2_H = \frac{1 - e^{-\alpha T} (-1)^j }{2} < 1 .
\]

A result corresponding to Theorem 3 for fractional Ornstein-Uhlenbeck sheets on $[0,T]^d$ with covariance structure
\begin{equation}
\E X_s X_t = \prod\limits^d_{i=1} e^{-\alpha_i |s_i - t_i|^{\rho_i}} ,\;
\alpha_i > 0, \;\rho_i \in (0,1]
\end{equation}
follows from Corollary 4.

Unfortunately, in the nonconvex case $\rho \in (1,2)$ it is not true that $\beta_j ( \rho) \geq 0$ for every $j \geq 0$ so that the approach of Theorem 3 does not work. In fact, starting again from
\[
\beta_j ( \rho) = \frac{2 \pi}{T} p_\rho ( \pi_j/T) - \frac{2}{T} \int^\infty_T \gamma^\rho (t) \cos ( \pi j t /T) dt,
\]
three integrations by parts show that
\[
\begin{array}{lcl}
\beta_j ( \rho) & = & \frac{2 \pi}{T} p_\rho ( \pi j/T) +
\frac{ (-1)^{j+1} 2 T \alpha \rho e^{- \alpha T^\rho} T^{\rho - 1} }{ (\pi j)^2 } + O ( j^{-3} ) \\ \\
& = & \frac{ (-1)^{j+1} 2 T \alpha \rho e^{- \alpha T^\rho} T^{\rho - 1} }{ (\pi j)2 } + O (
j^{-(1+\rho)} ) , j \rightarrow \infty.
\end{array}
\]
This means that for any $T > 0$, the $2T$-periodic extension of $\gamma ^\rho_{|[-T,T]}$ is not nonnegative definite for $\rho \!\in (1,2)$ in contrast to the case $\rho \! \in (0,1]$. \\

It is interesting to observe that the convex function $\gamma (t) = e^{- \alpha t^\rho} , \rho \in (0,1]$ in Theorem 3 can be replaced by any integrable convex positive function $\gamma$ on $(0,\infty)$. Let
$X = (X)_{t \in \R}$ be a continuous stationary centered Gaussian process with
$\E X_s X_t = \gamma (s-t)$, $\gamma : \R \rightarrow \R$. Then $\gamma$
is continuous, symmetric and nonnegative definite. Let
\[
\beta_0 := \frac{1}{2 T} \int^T_{-T} \gamma (t) dt, \; \beta_j := \frac{1}{T} \int^T_{-T} \gamma (t) \cos ( \pi j t/T) dt, \; j \geq 1 .
\]
Then the extension of Theorem 3 reads as follows.
\begin{Thm}
Assume that $\gamma$ is convex and positive on $(0, \infty)$ 
%and twice differentiable on $(0, \infty)$ such that
with
%$\gamma ( \infty) = 0$ and
$\gamma \in L1 ([T, \infty), dt)$.
Then $X$ admits a spectral density $p$ and we assume that $p$ satisfies $p \in L^2 (\R, dx)$ and
the high-frequency condition
\[
p(x) \sim c x^{-\delta} \; \mbox{as} \; x \rightarrow \infty
\]
for some $\delta \in (1,2], c \in (0, \infty)$. Then $\beta_j \geq 0$ for every $j \geq 0$ and the sequence
\begin{eqnarray}
f_0 & = & \sqrt{\beta_0} ,\quad  f_{2j} (t) = \sqrt{\beta_j} \cos ( \pi j t/T), \\
f_{2 j-1} (t) & = & \sqrt{\beta_j} \sin ( \pi j t/T) ,\quad  j \geq 1 \nonumber
\end{eqnarray}
is admissible for the contiuous process $X$  in ${\cal C}([0,T])$.
Moreover,
\[
l_n (X) \approx n^{-(\delta - 1)/2} ( \log n)^{1/2} \; \mbox{as} \; n \rightarrow \infty
\]
and the sequence (4.11) is rate optimal.
\end{Thm}
{\bf Proof.} The function $\gamma$ is integrable and convex over $[T,\infty)$, $\gamma(\infty) := \lim_{t\to \infty}\gamma(t)=0$ so that  
%Since
$\gamma$ is in fact a Polya-type function. Hence its right derivative $\gamma'$ is non-decreasing with $\gamma{'}(\infty) = 0$, the spectral measure of $X$ admits a Lebesgue-density $p$  and
\[
\gamma (t) = \int \left(1 - \frac{|t|}{s}\right)^{+} d \nu (s)
\]
for all $t \in \R$, where $\nu$ is a finite Borel measure on $(0, \infty)$ (see [10], Theorems 4.3.1 and 4.3.3).
Therefore, using Fubini's theorem,
it is enough to show the positivity of the numbers $\beta_j$ for functions of the type $\gamma (t) = (1 - \frac{|t|}{s} )^{+} , s \in (0, \infty)$. But in this case an integration by parts yields
\[
\beta_0 = \frac{T \wedge  s}{T} ( 1 - \frac{T \wedge  s}{2 s}) \geq 0 \; \mbox{ and } \; \beta_j = \frac{2T}{s (\pi j)^2}
(1 - \cos ( \pi j(T \wedge  s)/T)) \geq 0,\quad  j \geq 1.
\]
Now one proceeds along the lines of the proof of Theorem 3.
Since $\gamma$ is of bounded variation on $[-T,T]$, the representation (4.4) of $\E X_s X_t$ is true with $\beta_j(\rho)$ replaced by $\beta_j$ so that the sequence $(f_j)$ is admissible for $X$ in ${\cal C} ([0,1])$. Using $\gamma \in L1 (\R,dt)$ and the Fourier inversion formula, one gets for $j \geq 1$
\[
\beta_j = \frac{2 \pi}{T} p(\pi j /T) - \frac{2}{T} \int^\infty_T \gamma(t) \cos (\pi j t /T) dt.
\]
Since $\gamma ( \infty) = \gamma{'} (\infty) = 0$, integrating twice by parts yields
\[
\int^\infty_T \gamma(t) \cos ( \pi j t/T) dt = O(j^{-2}),
\]
hence
\[
\beta_j = O(j^{-\delta}) \; \mbox{as} \; j \rightarrow \infty
\]
in view of $\delta \leq 2$. Furthermore, the assumption $p \in L^2 ( \R, dx)$ and the high-frequency condition yield
\[
\lambda_j \sim c_1 j^{-\delta} \; \mbox{as} \; j \rightarrow \infty
\]
for an appropriate constant $c_1 \in (0, \infty)$ (see [16]).
Now, one derives from Proposition 4 the remaining assertions. \hfill{$\Box$}  \\ \\
{\bf Acknowlegdement.} It is a pleasure to thank Wolfgang Gawronski for helpful discussions on
Section 4. \\ \\
{\bf References}
\begin{itemize}
\item[{[1]}]
{\sc Ayache, A., Taqqu, M.S.}, Rate optimality of Wavelet series approximations of fractional Brownian motion,
{\em J. Fourier Anal. and Appl.},  {\bf 9} (2003), 451-471.
\item[{[2]}]
{\sc Bogachev, V.I.}, {\em Gaussian Measures}, AMS, 1998.
\item[{[3]}]
{\sc Bronski, J.C.}, Small ball constants and tight eigenvalue asymptotics for fractional Brownsian motions, {\em J. Theoret. Probab.} 16 (2003), 87-100.
\item[{[4]}]
{\sc Christensen, O.}, {\em An Introduction to Frames and Riesz Bases}, Birkh{\"a}user, Boston, 2003.
\item[{[5]}]
{\sc Dzhaparidze, K., van Zanten, H.}, A series expansion of fractional Brownian motion, {\em Probab. Theory
Relat. Fields},  {\bf 130} (2004), 39-55.
\item[{[6]}]
{\sc Dzhaparidze, K., van Zanten, H.}, Optimality of an explicit series expansion of the fractional Brownian sheet, {\em Statist. Probab. Letters} {\bf 71} (2005), 295 - 301.
\item[{[7]}]
{\sc Dzhaparidze, K., van Zanten, H.}, Krein's spectral theory and the Paley-Wiener expansion of fractional Brownian
motion, {\em Ann. Probab.} {\bf 33} (2005), 620-644.
\item[{[8]}]
{\sc K{\"u}hn, T., Linde, W.}, Optimal series representation of fractional Brownian sheets,
{\em Bernoulli} 8 (2002), 669 - 696.
\item[{[9]}]
{\sc Ledoux, M., Talagrand, M.}, {\em Probability in Banach Spaces}, Springer, Berlin, 1991.
\item[{[10]}]
{\sc Lukacs, E.}, {\em Characteristic Functions}, second edition, Griffin, London, 1970.
\item[{[11]}]
{\sc Luschgy, H.}, Linear estimators and radonifying operators, {\em Theory Probab. Appl.} {\bf 40} (1995), 205-213.
\item[{[12]}]
{\sc Luschgy, H., Pag{\`e}s, G.}, Sharp asymptotics of the functional quantization problem for Gaussian processes, {\em Ann. Probab.} 32 (2004), 1574- 1599.
\item[{[13]}]
{\sc Luschgy, H., Pag{\`e}s, G.}, High-resolution product quantization for Gaussian processes under sup-norm distortion, {\em Bernoulli} 13 (2007), 653-671.
\item[{[14]}]
{\sc Meyer, Y., Sellan, F. Taqqu, M.S.}, Wavelets, generalized white noise and fractional integration: the synthesis of fractional Brownian motion, {\em J. Fourier Anal. and Appl.} {\bf 5} (1999), 465-494.
\item[{[15]}]
{\sc Neveu, J.}, {\em Processus Al\'{e}atoires Gaussiens}, Les Presses de l'Universit\'{e} de
Montr\'{e}al, 1968.
\item[{[16]}]
{\sc Rosenblatt, M.,}  Some results on the asymptotic behaviour of eigenvalues for a class of integral equations with translation kernel, {\em  J. Math. Mech.} 12 (1963), 619 - 628.
\item[{[17]}]
{\sc Vakhania, N.N., Tarieladze, V.I., Chobanyan}, S.A., {\em Probability Distributions on Banach Spaces}, Kluwer, Boston, 1987.
\end{itemize}

\end{document}